\documentclass[12pt,a4paper]{amsart}
\usepackage{a4}
\usepackage[english]{babel} 
\usepackage{amsmath,amssymb,amsthm}
\def\leq {\leqslant}
\def\le {\leqslant}

\def\geq {\geqslant}
\def\@bibitem[#1]#2{\item\@biblabel{#1}.\if@filesw
{\def\protect##1{\string##1\space}\immediate\write
\@auxout{\string\bibcite{#2}{#1}}}\fi\ignorespaces\@showtag{#2}}

\theoremstyle{plain}
\newtheorem{theorem}{Theorem}
\newtheorem{rem}{Remark}
\newtheorem{lemma}{Lemma}
\renewcommand{\theequation}%
{\arabic{section}.\arabic{equation}}
\begin{document}

\title{Estimates of the order of approximation of functions of several variables in the generalized Lorentz space}
\author{ G. Akishev}
\address{ Lomonosov Moscow University, Kazakhstan Branch \\
Str. Kazhymukan, 11 \\
010010, Nur-Sultan, Kazakhstan}

\address{
Institute of mathematics and mathematical modeling\\
Pushkin str, 125 \\
050010, Almaty \\
 Republic of Kazakhstan
}

\maketitle

\begin{quote}
\noindent{\bf Abstract.}
In this paper we consider $ X(\bar\varphi)$ anisotropic symmetric space $ 2\pi$ of periodic functions of $m$ variables, in particular, the generalized Lorentz space $L_{\bar{\psi},\bar{\tau}}^{*}(\mathbb{T}^{m})$ and Nikol'skii--Besov's class $S_{X(\bar{\varphi}),\bar{\theta}}^{\bar r}B$.  
The article proves an embedding theorem for the Nikol'skii - Besov class in the generalized Lorentz space and establishes an upper bound for the best approximations by trigonometric polynomials with harmonic numbers from the hyperbolic cross of functions from the class $S_{X(\bar{\varphi}),\bar{\theta}}^{\bar r}B$.
\end{quote}
\vspace*{0.2 cm}



\section{Introduction}
\label{sec1}

Let $\mathbb{R}^{m}$ --- $m$ -- dimensional Euclidean point space $\bar{x} = (x_{1}, \ldots, x_{m})$ with real coordinates; $I^{m} = \{\bar{x} \in \mathbb{R}^{m}; \ 0 \leq x_{j} \leq 1; \ j = 1, \ldots, m \}=[0, 1]^{m} $ --- $ m $ -- dimensional
cub.
\smallskip

Two nonnegative Lebesgue measurable functions ~ $ f, $ $ g $ are called equimeasurable if
$$
\mu\{\bar x\in I^{m}\colon f(\bar x) > \lambda\} = \mu\{\bar x\in I^{m}\colon g(\bar x) > \lambda\},\quad \lambda > 0,
$$
where $\mu e $ is the Lebesgue measure of the set $e \subset I^{m}$.

For a nonnegative measurable function $ f $, a nonincreasing rearrangement is a function $f^{*}(t)=\inf\{\lambda >0 : \,\, \mu\{\bar x\in I^{m}\colon f(\bar x) > \lambda\}< t\}$. It is known that the functions $f$, $f^{*}$ are equimeasurable ([1], Ch. 2, Section 2).

Let $X$ be the Banach space of Lebesgue measurable functions $f$ on $I^{m}$ with norm $\|f\|_{X}$. The space $X$ is called symmetric if

1) from the fact that $|f(\bar x)|\leq |g(\bar x)|$ almost everywhere on $I^{m}$
and $g \in X $, it follows that $ f \in X $ and $ \| f \|_{X} \leq \| g \|_{X}$

2) from the fact that $ f \in X $ and the equimeasurability of the functions $|f(\bar x)|$ and $|g(\bar x)|$
it follows that $ g \in X $ and $ \|f\|_{X} = \|g\|_{X}$ (see\cite[Ch. 2, Sec.  4.1]{1}).

The norm $\|\chi_{e}\|_{X}$ of the characteristic function $\chi_{e}(t)$ of the measurable set $e \subset I^{m}$ is called
   is the fundamental function of the space $X$ and is denoted by $\varphi(\mu e) = \|\chi_{e}\|_{X} $.

It is known that the non-increasing rearrangement of the characteristic function $ \chi_{e}$ of the measurable set $ e \subset [0, 1]^{m}$ is equal to the function $ \chi_{[0, t]}$, where $ t = \mu e$. Therefore, the fundamental function of the symmetric space $ X $ is the function $\varphi(t) = \|\chi_{[0, t]}\|_{X}$, defined on the segment $[0, 1]$.
She  is a concave, non-decreasing, continuous function on $[0, 1]$, and
$\varphi (0) = 0$ (see \cite[Ch. 2, Sec.  4.4]{1}). Such functions are called $ \Phi$ -- functions.

For a given function $\varphi(t),$\ $t\in [0,1]$, we define 
$$
\alpha_{\varphi}={\underline\lim}_{t\rightarrow
0}\frac{\varphi(2t)}{\varphi(t)},\quad
\beta_{\varphi}=\overline{\lim}_{t\rightarrow 0}\frac{\varphi(2t)}{\varphi(t)}.
$$
Symmetric space $X$ with fundamental function $\varphi$ and norm
$\|f\|_{X}$ will be denoted by $X(\varphi) $ and its norm as $\|f\|_{X(\varphi)} $.

It is known that for any symmetric space $X(\varphi)$
inequalities  $1 \leq \alpha_{\varphi} \leq \beta_{\varphi} \leq 2$.

One example of a symmetric space is $ L_{q}(\mathbb{T}^{m})$--the Lebesgue space with the norm
$$
\|f\|_{q}=\biggl(\,\int\limits_{I^{m}}|f(2\pi\bar{x})|^{q}d\bar{
x}\,\biggr)^{1/q},\quad 1\le q < \infty.
$$
Here and in what follows, $\mathbb{T}^{m} = [0, 2\pi]^{m}$, the functions $f$ are $2\pi$--periodic in each variable.

\smallskip

Let the function  $\psi$ be continuous, concave and non-decreasing by $[0, 1]$, $\psi(0)=0$ and $0< \tau <\infty$.  
The generalized Lorentz space $L_{\psi, \tau}(\mathbb{T}^{m})$ is the set of measurable functions $f(\overline{x})=f(x_1,\ldots, x_m)$ of periodic $2\pi$ in each variable, such that   (see \cite{3})
$$
\|f\|_{\psi,\tau}^{*}=\bigg(\int\limits_{0}^{1} f^{*^\tau}(t) \psi^{\tau}(t)\frac{dt}{t}\bigg)^{1/\tau} <\infty.
 $$

It is known that under the conditions $ 1 <\alpha_{\psi}, \beta_{\psi} <2$, the space $ L_{\psi, \tau}(\mathbb{T}^{m})$ is a symmetric space with a fundamental function $\psi$.

Note that for  $\psi (t) = t^{1/q}$ the space $L_{\psi, \tau}(\mathbb{T}^{m})$ coincides with the Lorentz space  $L_{q, \tau}(\mathbb{T}^{m})$, $1<q, \tau <\infty$, which consists of all functions $f$ such that (see \cite[Ch. 1, Sec. 3]{2}) 
$$
\|f\|_{q,\tau}=\Bigg(\frac{\tau}{q}\int\limits_{0}^{1}\biggl(\int\limits_{0}^{t}f^{*}(y)dy\biggr)^{\tau}t^{\tau(\frac{1}{q}-1)-1}dt
\Bigg)^{1/\tau} < \infty.
$$

We consider $ X(\bar\varphi)$ anisotropic symmetric space $ 2\pi$ of periodic functions of $ m $ variables, with the norm $\|f\|_{X(\bar\varphi)}^{*} = \|\ldots\|f^{*_{1},...,*_{m}}\|_{X(\varphi_{1})}\ldots\|_{X(\varphi_{m})}$, where
$f^{*_{1},...,*_{m}}(t_{1},...,t_{m})$
non-increasing rearrangement  of a function
$|f(2\pi \bar{x})|$ for each variable $x_{j} \in [0, 1]$ with fixed other variables (see [5]) and $X(\varphi_{j})$ --- symmetric space in the variable $ x_{j}$, with the fundamental function $\varphi_{j}$ (see [4]).

The associated space to the symmetric space $X(\bar\varphi)$ is the space of all measurable functions $ g $ for which  (see  \cite{4})
$$
\sup_{{}_{\|f\|_{X(\bar\varphi)}^{*}\leq 1}^{f\in X(\bar\varphi)}} \int_{I^{m}}
f(2\pi\bar{x})g(2\pi\bar{x})d\bar{x} < \infty, 
 $$
and is denoted by the symbol $X^{'}(\bar{\tilde{\varphi}})$, and its norm as $\|g\|_{X^{'}(\bar{\tilde{\varphi}})}^{*}$, 
where  
$$
\bar{\tilde{\varphi}}(t)=(\tilde{\varphi}_{1}(t),...,\tilde{\varphi}_{m}(t)),  \,\, \tilde{\varphi}_{j}(t)=\frac{t}{\varphi_{j}(t)}
$$
 for $t\in (0, 1]$ and $\tilde{\varphi}_{j}(0)=0$, $j=1,...,m$.
It is known that 
$$
\left|\int_{I^{m}}
f(2\pi\bar{x})g(2\pi\bar{x})d\bar{x}\right| \leq \|g\|_{X^{'}(\bar{\tilde{\varphi}})}^{*} \|f\|_{X(\bar\varphi)}^{*}, \,\, f\in X(\bar\varphi), g\in X^{'}(\bar{\tilde{\varphi}}).
$$

Let $\bar{x}=(x_{1},...,x_{m})\in \mathbb{I}^{m}=[0, 1]^{m}$ and given 
$\Phi$--functions $\psi_{j}(x_{j})$, $x_{j}\in[0, 1]$ and
$\tau_{j}\in[1,+\infty)$, $j=1,...,m.$ We shall denote by $L_{\bar{\psi},\bar{\tau}}^{*}(\mathbb{T}^{m})$ the generalized Lorentz space with anisotropic norm of Lebesgue measurable functions $f(2\pi \bar{x})$ of period $2\pi$ in each variable , such that the quantity   
$$
\|f\|_{\bar{\psi},\bar{\tau}}^{*}=\Bigl[\int_{0}^{1}\psi_{m}^{\tau_{m}}
(t_{m})
\Bigl[\ldots\Bigl[\int_{0}^{1}\psi_{1}^{\tau_{1}}(t_{1})\left(
f^{*_{1},...,*_{m}}(t_{1},...,t_{m})
\right)^{\tau_{1}}\frac{dt_{1}}{t_{1}}\Bigr]^
{\frac{\tau_{2}}{\tau_{1}}}
\Bigr]^{\frac{\tau_{m}}{\tau_{m-1}}}\frac{dt_{m}}{t_{m}}
\Bigr]^{\frac{1}{\tau_{m}}}
$$
is finite.
For functions
$\psi_{j}(t)=t^{\frac{1}{q_{j}}}$, $j=1,...,m$ Lorentz space
$L_{\bar{\psi},\bar{\tau}}^{*}(\mathbb{T}^{m})$, will be denoted by
$L_{\bar{q},\bar{\tau}}^{*}(\mathbb{T}^{m})$ and instead of $\|\bullet\|_{\bar{\psi}, \bar{\tau}}^{*}$  respectively we will write
$\|\bullet\|_{\bar{q},\bar{\tau}}^{*}$ (see \cite{5}).
 
Let $\overset{\circ \;\;}
L_{\bar{\psi}, \bar{\tau}}^{*}
\left(\mathbb{T}^{m} \right)$ be the set of functions $f\in
L_{\bar{\psi}, \bar{\tau}}^{*}(\mathbb{T}^{m})$ such
that
 $$
\int\limits_{0}^{2\pi }f\left(\overline{x} \right) dx_{j}
=0,\;\;\forall j=1,...,m .
 $$
We will use the following notation:
$a_{\overline{n} } (f)$ be the Fourier coefficients of $f\in
L_{1}(\mathbb{T}^{m})$ with respect to the multiple
trigonometric system and 
$$
\delta _{\overline{s} } \left( f,\overline{x} \right)
=\sum\limits_{\overline{n} \in \rho \left( \overline{s} \right)
}a_{\overline{n} } \left( f\right) e^{i\langle\overline{n} ,\overline{x}\rangle } ,
$$
where $\langle\bar{y},\bar{x}\rangle=\sum\limits_{j=1}^{m}y_{j}
x_{j}$,
 $$
\rho (\bar{s})=\left\{ \overline{k} =\left( k_{1}
,...,k_{m} \right) \in \mathbb{Z}^{m}: \quad 2^{s_{j} -1} \leq \left| k_{j}
\right| <2^{s_{j} } ,j=1,...,m\right\}.
$$

We will consider the functional class of Nikol'skii-Besov
 $$
S_{X(\bar{\varphi}),\bar{\theta}}^{\bar r}B=
\Bigl\{f\in \overset{\circ \;\;}X(\bar{\varphi}) : \quad
\|f\|_{X(\bar{\varphi})}^{*} + \Bigl\|\Bigl\{\prod_{j=1}^{m}
2^{s_{j}r_{j}} \|\delta_{\bar
s}(f)\|_{X(\bar{\varphi})}^{*} \Bigr\}_{\bar{s}\in
\mathbb{Z}_{+}^{m}}\Bigr\|_{l_{\bar\theta}}\leq 1\Bigr\},
 $$
where $\bar{\theta}=(\theta_{1},...,\theta_{m}),$  $\bar{r}=(r_{1},...,r_{m}),$ 
$1\leq\theta_{j}\le+\infty,$ $0<r_{j}<+\infty,$  $j=1,...,m.$

In the case of $X(\bar{\varphi})=L_{p}(\mathbb{T}^{m})$, $1\leq p < \infty$, $ 1 \leq p <\infty $, class $S_{X(\bar{\varphi}),\bar{\theta}}^{\bar r}B$ is defined and studied in \cite{6}--\cite{8}.

For a fixed vector $\bar{\gamma}=(\gamma_{1},\ldots,\gamma_{m}),$ \;\;$\gamma_{j}>0, \;\; j=1,\ldots,m$, set 
   $$
Q_{n}^{\bar\gamma}=
 \cup_{{}_{\langle\bar{s},\bar{\gamma}\rangle < n}}\rho(\bar{s}), \quad
 T(Q_{n}^{\bar \gamma})=\{t(\bar x)=\sum\limits_{\bar{k}\in
Q_{n}^{\bar\gamma}}b_{\bar k} e^{i\langle\bar{k},\bar{x}\rangle}\},
 $$

$E_{n}^{(\overline{\gamma})}(f)_{X(\bar{\varphi})}$
is the best approximation of a
function $f\in X(\bar{\varphi})$ by polynomials in $T(Q_{n}^{\bar\gamma})$ , and
$S_{n}^{\bar \gamma}(f,\bar{x})=\sum_{\bar{k}\in Q_{n}^{\bar
\gamma}}a_{\bar k}(f)\cdot e^{i\langle\bar{k}, \bar{x}\rangle}$ is
a partial sum of the Fourier series of $f$.

We shall denote by  $C(p,q,y,..)$  positive quantities which depend
only on the parameter in the parentheses and not necessarily the
same in distinct formulae . The notation $A\left( y\right) \asymp
B\left( y\right)$ means that there exist positive constants
 $C_{1},\,C_{2} $ such that  $C_{1} \cdot A\left( y\right)
\leq B\left( y\right) \leq C_{2} \cdot A\left( y\right)$.

Exact order estimates for the best approximation of functions of various classes in the Lebesgue space $L_{p}(\mathbb{T}^{m})$ are well known (see survey articles \cite{9}--\cite{11}, monograph \cite{12} and bibliographies in them). These questions in the space $L_{\bar{q}, \bar{\tau}}^{*}(\mathbb{T}^{m})$ were studied in \cite{13}--\cite{18}.

The main aim of the present paper is to find the order of the quantity 
$$
E_{n}^{(\bar\gamma)}(S_{X(\bar{\varphi}),\bar{\theta}}^{\bar{r}}B)
_{\bar{\psi},\bar{\tau}}\equiv \sup_{f\in
S_{X(\bar{\varphi}),\bar{\theta}}^{\bar r}B}E_{n}^{(\overline{\gamma})}(f)_{\bar{\psi},\bar{\tau}}.
$$

This paper is organized as follows. In  Section 1 we give auxiliary results.  
 In Section 2, we will prove the main results. Our main results in this Section reads:

{\bf Theorem 2.} \textit{ Let $1<\alpha_{\psi_{j}}\leq \beta_{\psi_{j}}<
\alpha_{\varphi_{j}}\leq \beta_{\varphi_{j}}<2$, $1\leq\tau_{j}<+\infty$,
$j=1,...,m$. If $f\in X(\bar{\varphi})$ and
 $$
\Bigl\{\prod_{j=1}^{m}\frac{\psi_{j}(2^{-s_{j}})}{\varphi_{j}
(2^{-s_{j}})}\|\delta_{\bar s}(f)\|_{X(\bar{\varphi})}^{*}
\Bigr\}_{\bar{s}\in{\mathbb Z}_{+}^{m}}\in
l_{\bar{\tau}},
 $$
 then $f\in \overset{\circ}{L}_{\bar{\psi},\bar{\tau}}^{*}(\mathbb{T}^{m})$ and the following inequality holds 
 $$
\|f\|_{\bar{\psi},\bar{\tau}}^{*}\leq C
\Bigl\|\Bigl\{\prod_{j=1}^{m}\frac{\psi_{j}(2^{-s_{j}})}{\varphi_{j}
(2^{-s_{j}})}\|\delta_{\bar s}(f)\|_{X(\bar{\varphi})}^{*}
\Bigr\}_{\bar{s}\in{\mathbb Z}_{+}^{m}}
\Bigr\|_{l_{\bar\tau}}.
 $$
 }

Further, for brevity, we put
 $
\mu_{j}(s)=\frac{\psi_{j}\left(2^{-s}\right)}
{\varphi_{j}\left(2^{-s}\right)}.
 $

{\bf Theorem 5.}
{\it Let $1\leq\theta_{j}\leq +\infty,$ $1\leq\tau_{j}<+\infty,$ $r_{j}>0$, $\frac{r_{j}}{r_{1}}$, $j=1,...,m$ and functions $\varphi_{j},$ $\psi_{j}$ satisfy the conditions  $1<\alpha_{\psi_{j}}
\leq\beta_{\psi_{j}}<\alpha_{\varphi_{j}}\leq\beta_{\varphi_{j}}<2,$ $j=1,...,m$
and 
$$
\Bigl[\sum_{s_{j}=0}^{\infty}\Bigl(\frac{\psi_{j}(2^{-s_{j}})}{\varphi_{j}(2^{-s_{j}}
)}2^{-s_{j}r_{j}}\Bigr)^{\varepsilon_{j}}\Bigr]^{\frac{1}{\varepsilon_{j}}}<+\infty,
$$
where $\varepsilon_{j}=\tau_{j}\beta_{j}'$, $\beta_{j}'=
\frac{\beta_{j}}{\beta_{j}-1}$, $j=1,...,m$, if $\beta_{j}=
\frac{\theta_{j}}{\tau_{j}}>1$ and $\varepsilon_{j}=+\infty$, if
$\theta_{j}\leq\tau_{j}<\infty$, $j=1,...,m.$

1) If $1\leq\tau_{j}<\theta_{j}<+\infty$, $j=1,...,m,$ then
 $$
E_{n}^{(\bar\gamma)}(S_{X(\bar{\varphi}),\bar{\theta}}^{\bar{r}}B)
_{\bar{\psi},\bar{\tau}}\leq
C\Bigl\|\Bigl\{\prod_{j=1}^{m}2^{-s_{j}r_{j}}
\mu_{j}(s_{j})\Bigr\}_{\bar{s}\in Y^{m}(\bar{\gamma},n)}
\Bigr\|_{l_{\bar\varepsilon}}.
 $$

2) If $1\leq \theta_{j}\leq\tau_{j}<+\infty,$  $j=1,...,m,$ then
 $$
E_{n}^{\bar \gamma}(S_{X(\bar{\varphi}),\bar{\theta}}^{\bar r}B)
_{\bar{\psi},\bar{\tau}}\leq
C\sup\Bigl\{\prod_{j=1}^{m}2^{-s_{j}r_{j}}
\mu_{j}(s_{j}): \quad \bar{s}\in\mathbb{Z}_{+}^{m},
 \langle\bar{s},\bar{\gamma}\rangle\geq n\Bigr\}.
 $$}
Here $Y^{m}(\bar{\gamma},n)=\{\bar{s}\in\mathbb{Z}_{+}^{m},
 \langle\bar{s},\bar{\gamma}\rangle\geq n \}$. 

\section{Auxiliary statements}
\label{sec 1}

In this section, we present some well-known results and prove several lemmas.
The following generalization of the discrete Hardy inequality is known
 
\begin{lemma}\label{lem1} (see \cite{19}).  
{\it Let $0<\theta<+\infty$ and given positive numbers
 $a_{k}$, $b_{k}, k=0,1,2,...$ .

a) If
$
\sum\limits_{k=0}^{n}a_{k}\leq C\cdot a_{n},\,\,
$
then
$
\sum\limits_{n=0}^{\infty}a_{n}\Bigl(\sum\limits_{k=n}^{\infty}
b_{k}\Bigr)^{\theta}\leq
C\cdot\sum\limits_{n=0}^{\infty}a_{n}b_{n}^{\theta}.
$

b) If
$
\sum\limits_{k=n}^{\infty}a_{k}\leq C a_{n},\,\,
$
then
$
\sum\limits_{n=0}^{\infty}a_{n}\Bigl(
\sum\limits_{k=0}^{n}b_{k}\Bigr)^{\theta}\leq C\cdot
\sum\limits_{n=0}^{\infty}a_{n}\cdot b_{n}^{\theta}.
$
}
\end{lemma}

The proof of Lemma 1 is also given in \cite{24}.

 \begin{lemma}\label{lem2} (see \cite{20}, \cite{21}) 
 {\it If a
$1<\alpha_{\psi},\beta_{\psi}<2$ for $\Phi$--function
$\psi (x)$, $x \ in [0,1] $, then for any $ q> 0 $
  the relations are satisfied
   $$ \int\limits_{0}^{x}\frac{\psi^{q}(t)}{t}dt
=O(\psi^{q}(x)), x\rightarrow +0, 
$$
$$
\int\limits_{x}^{1}[t\psi^{q}(t)]^{-1}dt=O(\psi^{-q}(x)),
x\rightarrow +0.
$$
}
\end{lemma}

\begin{lemma}\label{lem3} (see \cite{20}, \cite{21}) 
{\it If $\Phi$-th  functions
$\varphi(x), \psi(x),x\in[0,1]$ satisfies the condition
$\alpha_{\varphi}>\beta_{\psi}>1$, then for
function 
$$
g(x)=\left\{
\begin{array}{ll}\frac{\varphi(x)}{\psi(x)}, &\mbox{if} \,\, x\in(0,1]\\
0, &\mbox{if} \,\,  x=0
\end{array}\right.
$$
there is a $\Phi$-th function $g_{1}(x)$
such that
 $g(x)\asymp g_{1}(x),x\in [0,1]$ and
$\alpha_{g_{1}}>1$.
}
\end{lemma}

\begin{lemma}\label{lem4}  
{\it If $\Phi$-th  functions
$\varphi(x), \psi(x),x\in[0,1]$ satisfies the condition 
$1<\alpha_{\psi} \leq \beta_{\psi} <\alpha_{\varphi}\leq \beta_{\varphi}<2$ and $0< \theta < \infty$, then the following inequality holds 
 $$
\sum\limits_{s=0}^{n}\Bigl(\frac{\psi(2^{-s})}{\varphi(2^{-s})}\Bigr)^{\theta} \leq C \Bigl(\frac{\psi(2^{-n})}{\varphi(2^{-n})}\Bigr)^{\theta}, \,\, n\in \mathbb{N}.
$$
}
\end{lemma}
{\bf Proof.} We will consider the function  
$$
g(x) = \left\{
\begin{array}{ll}\frac{\varphi(x)}{\psi(x)}, &\mbox{if} \,\, x\in(0,1]\\
0 , &\mbox{if} \,\,  x=0 .
\end{array}\right.
$$
Since $1<\alpha_{\psi} \leq \beta_{\psi} <\alpha_{\varphi}$, then according to Lemma 3 there exists a $\Phi$ - function $ g_{1}$ such that $g(x) \asymp g_{1}(x) $, $ x \rightarrow + 0 $ and $\alpha_{g_{1}}> 1$.
Therefore  
 $$
\sum\limits_{s=0}^{n}\Bigl(\frac{\psi(2^{-s})}{\varphi(2^{-s})}\Bigr)^{\theta} \leq C\sum\limits_{s=0}^{n}g_{1}^{-\theta}(2^{-s}). \eqno (1)
$$
Since $g_{1}\uparrow$ on $(0, 1]$, then
$$
\sum\limits_{s=0}^{n}g_{1}^{-\theta}(2^{-s})\leq \ln 2 \sum\limits_{s=0}^{n}\int\limits_{2^{-s-1}}^{2^{-s}} g_{1}^{-\theta}(t)\frac{dt}{t} = \ln 2 \int\limits_{2^{-n-1}}^{1} g_{1}^{-\theta}(t)\frac{dt}{t}  \eqno (2)
$$
Now, from inequalities (1) and (2), according to Lemma 2, we obtain the assertion of the lemma.

\begin{lemma}\label{lem5} 
{\it Let $\Phi$ be a function
$\psi$ satisfies the conditions $1<\alpha_{\psi} \leq \beta_{\psi} <2$, then the following inequality holds
 $$
\sum\limits_{s=n}^{\infty}\psi^{\theta}(2^{-s}) \leq C\psi^{\theta}(2^{-n}), n\in \mathbb{N}. 
$$
}
\end{lemma}
{\bf Proof.} Since $\frac{\psi(t)}{t} \downarrow$ on $(0, 1]$, then
$$
\psi^{\theta}(2^{-s}) \leq C\int\limits_{2^{-s-1}}^{2^{-s}} \psi^{\theta}(t)\frac{dt}{t}.
$$
Therefore 
$$
\sum\limits_{s=n}^{\infty}\psi^{\theta}(2^{-s}) \leq C\sum\limits_{s=n}^{\infty}
\int\limits_{2^{-s-1}}^{2^{-s}} \psi^{\theta}(t)\frac{dt}{t} = C\int\limits_{0}^{2^{-n}} \psi^{\theta}(t)\frac{dt}{t}.
$$
Hence, according to Lemma 2, from this we obtain the assertion of the lemma.

{\bf Remark.} In Lemma 4 and Lemma 5 the coefficients C do not depend on $ n $. These lemmas were proved in integral form in \cite{20}.

We now prove a multidimensional version of Lemma 1.

\begin{lemma}\label{lem6} 
{\it Let positive numbers $b_{\bar{k}}=b_{k_{1},\ldots , k_{m}}$ be given for $\bar{k}=(k_{1},\ldots , k_{m})\in \mathbb{Z}_{+}^{m}$ and
  $a_{k_{j}}$, $ k_{j}=0,1,2,...$ and $1 \leq \theta_{j}<+\infty$, $j=1,\ldots , m$.
 
a) If
$$
\sum\limits_{k_{j}=0}^{n_{j}}a_{k_{j}}\leq C a_{n_{j}},\,\, j=1,\ldots , m,
$$
then
$$
A_{\theta_{1},\ldots ,\theta_{m}}=\left\{\sum\limits_{n_{m}=0}^{\infty}a_{n_{m}}\left[\sum\limits_{n_{m-1}=0}^{\infty}a_{n_{m-1}} \ldots \left[\sum\limits_{n_{1}=0}^{\infty}a_{n_{1}}\left( \sum\limits_{k_{m}=n_{m}}^{\infty}\ldots\sum\limits_{k_{1}=n_{1}}^{\infty}b_{\bar{k}} \right)^{\theta_{1}}\right]^{\frac{\theta_{2}}{\theta_{1}}}\ldots \right]^{\frac{\theta_{m}}{\theta_{m-1}}}\right\}^{\frac{1}{\theta_{m}}}
$$
$$
\leq C \left\{\sum\limits_{n_{m}=0}^{\infty}a_{n_{m}}\left[\sum\limits_{n_{m-1}=0}^{\infty}a_{n_{m-1}} \ldots \left[\sum\limits_{n_{1}=0}^{\infty}a_{n_{1}}b_{\bar{n}}^{\theta_{1}}\right]^{\frac{\theta_{2}}{\theta_{1}}}\ldots \right]^{\frac{\theta_{m}}{\theta_{m-1}}}\right\}^{\frac{1}{\theta_{m}}}.
$$

b) If 
$$
\sum\limits_{k_{j}=n_{j}}^{\infty}a_{k_{j}}\leq Ca_{n_{j}},\,\,j=1,\ldots , m,
$$
then
$$
\left\{\sum\limits_{n_{m}=0}^{\infty}a_{n_{m}}\left[\sum\limits_{n_{m-1}=0}^{\infty}a_{n_{m-1}} \ldots \left[\sum\limits_{n_{1}=0}^{\infty}a_{n_{1}}\left( \sum\limits_{k_{m}=0}^{n_{m}}\ldots\sum\limits_{k_{1}=0}^{n_{1}}b_{\bar{k}} \right)^{\theta_{1}}\right]^{\frac{\theta_{2}}{\theta_{1}}}\ldots \right]^{\frac{\theta_{m}}{\theta_{m-1}}}\right\}^{\frac{1}{\theta_{m}}}
$$
$$
\leq C \left\{\sum\limits_{n_{m}=0}^{\infty}a_{n_{m}}\left[\sum\limits_{n_{m-1}=0}^{\infty}a_{n_{m-1}} \ldots \left[\sum\limits_{n_{1}=0}^{\infty}a_{n_{1}}b_{\bar{n}}^{\theta_{1}}\right]^{\frac{\theta_{2}}{\theta_{1}}}\ldots \right]^{\frac{\theta_{m}}{\theta_{m-1}}}\right\}^{\frac{1}{\theta_{m}}}.
$$
}
\end{lemma}
{\bf Proof.} Let us prove item a). For $ m = 1 $ the statement is known (see Lemma 1). Let  $m=2$. Since $ 1 \leq \theta_{1} <\infty$, then by the property of the norm we have
$$
\sum\limits_{n_{2}=0}^{\infty}a_{n_{2}} \left[\sum\limits_{n_{1}=0}^{\infty}a_{n_{1}}\left( \sum\limits_{k_{2}=n_{2}}^{\infty}\sum\limits_{k_{1}=n_{1}}^{\infty}b_{\bar{k}} \right)^{\theta_{1}}\right]^{\frac{\theta_{2}}{\theta_{1}}} \leq \sum\limits_{n_{2}=0}^{\infty}a_{n_{2}} \left[\sum\limits_{k_{2}=n_{2}}^{\infty}\left(\sum\limits_{n_{1}=0}^{\infty}a_{n_{1}}\left(\sum\limits_{k_{1}=n_{1}}^{\infty}b_{\bar{k}} \right)^{\theta_{1}}\right)^{\frac{1}{\theta_{1}}}\right]^{\theta_{2}}
$$
Now, applying statement a) of Lemma 1 twice, from this we obtain
$$
\sum\limits_{n_{2}=0}^{\infty}a_{n_{2}} \left[\sum\limits_{n_{1}=0}^{\infty}a_{n_{1}}\left( \sum\limits_{k_{2}=n_{2}}^{\infty}\sum\limits_{k_{1}=n_{1}}^{\infty}b_{\bar{k}} \right)^{\theta_{1}}\right]^{\frac{\theta_{2}}{\theta_{1}}} \leq \sum\limits_{n_{2}=0}^{\infty}a_{n_{2}} \left[\sum\limits_{n_{1}=0}^{\infty}a_{n_{1}}b_{\bar{n}}^{\theta_{1}}\right]^{\frac{\theta_{2}}{\theta_{1}}}.
$$
Now suppose the statement is true for $m-1$. Let us prove it for $m$.
Since $\theta_{j}$, $ j = 1, \ldots, m-1$, then by the property of the norm and using the assumption, we have
$$
A_{\theta_{1},\ldots ,\theta_{m}}^{\theta_{m}} \leq 
$$
$$
\sum\limits_{n_{m}=0}^{\infty}a_{n_{m}}\left[\sum\limits_{k_{m}=n_{m}}^{\infty} \left[\sum\limits_{n_{m-1}=0}^{\infty}a_{n_{m-1}} \ldots \left[\sum\limits_{n_{1}=0}^{\infty}a_{n_{1}}\left( \sum\limits_{k_{m-1}=n_{m-1}}^{\infty}\ldots\sum\limits_{k_{1}=n_{1}}^{\infty}b_{\bar{k}} \right)^{\theta_{1}}\right]^{\frac{\theta_{2}}{\theta_{1}}}\ldots \right]^{\frac{1}{\theta_{m-1}}}\right]^{\theta_{m}} 
$$
$$
\leq C\sum\limits_{n_{m}=0}^{\infty}a_{n_{m}}\left[\sum\limits_{k_{m}=n_{m}}^{\infty} \left[\sum\limits_{n_{m-1}=0}^{\infty}a_{n_{m-1}} \ldots \left[\sum\limits_{n_{1}=0}^{\infty}a_{n_{1}}b_{n_{1},\ldots, n_{m-1}, k_{m}} ^{\theta_{1}}\right]^{\frac{\theta_{2}}{\theta_{1}}}\ldots \right]^{\frac{1}{\theta_{m-1}}}\right]^{\theta_{m}} 
$$
Now, by applying assertion a) of Lemma 1, from this we obtain the required inequality.
Assertion $ b) $ is proved similarly.
   
\begin{lemma}\label{lem7} {\it Let $1 < \alpha_{\varphi_{j}}\leq \beta_{\varphi_{j}} \leq 2,
$ $j=1,...,m$ and $E_{j}\subset
[0, 2\pi), j=1,...,m,$ be measurable sets. Then for any trigonometric polynomial
 $$
T_{\bar n}(\bar x)=\sum_{k_{1}=-n_{1}}^{n_{1}}...
\sum_{k_{m}=-n_{m}}^{n_{m}}b_{\bar k}e^{i\langle\bar{x},\bar{k}\rangle}  
 $$
the following inequality holds: 
 $$
\int_{E_{m}}...\int_{E_{1}} |T_{\bar n}(\bar x)|dx_{1}...dx_{m}\leq
 $$
 $$
\leq C\prod_{j\in e}
\frac{1}{\varphi_{j}(n_{j}^{-1})}\prod_{j\in e}|E_{j}|\prod_{j\in
\bar e} \bar{\varphi_{j}}(|E_{j}|)  
 \|T_{\bar n}\|_{X(\bar\varphi)}^{*},
 $$
where  $e\subset \{1,...,m\}$ and $\bar{e}$--complement set $e$ ,  $|E_{j}|$ --- Lebesgue measure of $E_{j}$.
}
\end{lemma}

{\bf Proof.} Let $e$ be an arbitrary subset of 
 $\{1,...,m\}.$ Also let  $E_{e}=\prod_{j\in e}E_{j},$
$d^{e}\bar x=\prod_{j\in e}dx_{j}.$
 It is known that for any trigonometric polynomial  $T_{\bar
n}(\bar x)$ for fixed  $x_{j}, j\in \bar e,$ the following formula holds:
 $$
T_{\bar n}(\bar x)=\frac{1}{(2\pi)^{|e|}}\int_{I^{|e|}} T_{\bar
n}(\bar y^{e},\bar x^{\bar e})\cdot \prod_{j\in
e}D_{n_{j}}(x_{j}-y_{j})d^{e}\bar y,
 $$
where  $D_{n}(t)$ -- is the Dirichlet kernel of the trigonometric system, 
$\bar{y}^{e}$ -- the vector with the coordinates $y_{j}$ if $j\in e$, and 
$\bar{x}^{\bar e}$ -- the vector with the coordinates $x_{j}$ if $j\in \bar
e$ ($\bar e$ is the set--theoretic completion of $e$ ) and $|e|$ is the number of elements of $e.$
By the property of Lebesgue integral, we have
 $$
\int_{E_{\bar e}}|T_{\bar n}(\bar x)|d^{\bar e}\bar x =
\frac{1}{(2\pi)^{|e|}}\int_{E_{\bar e}}|\int_{I^{|e|}} T_{\bar
n}(\bar y^{e},\bar x^{\bar e})\cdot \prod_{j\in
e}D_{n_{j}}(x_{j}-y_{j})d^{e}\bar y|d^{\bar e}\bar x\leq
 $$
 $$
\leq \frac{1}{(2\pi)^{|e|}}\int_{E_{\bar e}}\int_{I^{|e|}}
|T_{\bar n}(\bar y^{e},\bar x^{\bar e})|\cdot \prod_{j\in
e}|D_{n_{j}}(x_{j}-y_{j})|d^{e}\bar y d^{\bar e}\bar x =
 $$
 $$
=\frac{1}{(2\pi)^{|e|}}\int_{I^{m}} |T_{\bar n}(\bar y^{e},\bar
x^{\bar e})|\cdot \prod_{j\in e}|D_{n_{j}}(x_{j}-y_{j})| \cdot
\prod_{j\in \bar e}\chi_{E_{j}}(x_{j}) d^{e}\bar y d^{\bar e}\bar
x,
 $$
where $\chi_{E}$ is the characteristic function of the set $E.$
By appling the H\"{o}lder's inequality to the integrals in the right side of the last relation. Then 
 $$
\int_{E_{\bar e}}|T_{\bar n}(\bar x)|d^{\bar e}\bar x \leq
\frac{1}{(2\pi)^{|e|}} \|T_{\bar n}\|_{X(\bar\varphi)}^{*}
 \|\prod_{j\in e}|D_{n_{j}}(x_{j}-\bullet)|
\prod_{j\in \bar e}\chi_{E_{j}}\|_{X(\tilde{\bar\varphi})}^{*}=
 $$
 $$
=\frac{1}{(2\pi)^{|e|}} \|T_{\bar n}\|_{X(\bar\varphi)}^{*}
 \prod_{j\in e}\|D_{n_{j}}(x_{j}-\bullet)\|_{X(\tilde{\varphi}_{j})}^{*}
\prod_{j\in \bar
e}\|\chi_{E_{j}}\|_{X(\tilde{\varphi}_{j})}^{*}. \eqno(3)
 $$
 Since
 $$
\|\chi_{E_{j}}\|_{X(\tilde{\varphi}_{j})}^{*}=\tilde{\varphi}_{j}(|E_{j}|)
 $$
then from inequality (3) we obtain
$$
\int_{E_{\bar e}}|T_{\bar n}(\bar x)|d^{\bar e}\bar x \leq
\frac{1}{(2\pi)^{|e|}} \|T_{\bar n}\|_{X(\bar\varphi)}^{*}
$$
$$
\prod_{j\in \bar e}\tilde{\varphi}_{j}(|E_{j}|)  \prod_{j\in e} \sup_{x_{j}\in [0, 2\pi]}
\|D_{n_{j}}(x_{j}-\bullet)\|_{\tilde{\varphi}_{j}}^{*}, \eqno(4)
$$
where $d^{\bar e}\bar x = \prod\limits_{j\in e}dx_{j}$.

In the one-dimensional case, the estimate is known (see, for example, \cite{23})
 $$
\sup_{x_{j}\in [0, 2\pi]} \|D_{n_{j}}(x_{j}-\bullet)\|_{X(\tilde{\varphi}_{j})}^{*}\leq
C n_{j}\tilde{\varphi}_{j}(n_{j}^{-1}), \,\, j=1,\ldots m .  \eqno(5)
$$
Now from inequalities (4) and (5) we will have
$$
\int_{E_{\bar e}}|T_{\bar n}(\bar x)|d^{\bar e}\bar x \leq
C \|T_{\bar n}\|_{X(\bar\varphi)}^{*}
\prod_{j\in \bar e}\tilde{\varphi}_{j}(|E_{j}|) \prod_{j\in
e} n_{j}\tilde{\varphi}_{j}(n_{j}^{-1})  
$$
for fixed $x_{j}, j\in e.$
By integrating by the variables  $x_{j}, j\in e,$ the both sides of this inequality, we get 
  $$
\int_{E_{\bar e}}\int_{E_{e}}|T_{\bar n}(\bar x)|d^{e}\bar x
d^{\bar e}\bar x \leq C \|T_{\bar n}\|_{X(\bar\varphi)}^{*} \prod_{j\in \bar e}\tilde{\varphi}_{j}(|E_{j}|)\prod_{j\in
e} n_{j}\tilde{\varphi}_{j}(n_{j}^{-1})\prod_{j\in e}|E_{j}|.
$$
The proof is finished.

\section{Main results}\label{sec 2}

We now prove the main results.
We set
 $$
G_{e}(\bar n)=\{\bar{s}=(s_{1},...,s_{m})\in \mathbb{N}^{m}: s_{j}
\leq n_{j}, j\in e; \quad s_{j}>n_{j}, j\notin e\},
 $$
 where $e\subset\{1,...,m\}$;
 $$
U_{\bar n}(f,\bar{x})=\sum_{e\subset\{1,...,m\}}\sum_{\bar{s}\in G_{e}(\bar
n)} \delta_{\bar s}(f,\bar{x}).
 $$
Let
 $$
\bar{f}(\bar t)=\sup_{|E_{m}|\geq t_{m}}\frac{1}{|E_{m}|}
\int_{E_{m}}...\sup_{|E_{1}|\geq t_{1}}\frac{1}{|E_{1}|}
\int_{E_{1}}|f(x_{1},...,x_{m})|dx_{1}...dx_{m},
 $$
$|E_{j}|$ is the Lebesgue measure of the set $E_{j}\subset [0,2\pi).$

\begin{theorem}\label{th1} 
{\it Let  $\bar{\varphi}=(\varphi_{1},...,\varphi_{m})$ and  functions $\varphi_{j}$  satisfies the conditions $1<\alpha_{\varphi_{j}}\leq
\beta_{\varphi_{j}}<2$, $j=1,\ldots , m$. Then for each
function $ f\in X(\bar{\varphi})$ the next
inequality holds
 $$
{\bar{f}}(\bar{t})\leq
C\left\{\prod\limits_{j=1}^{m}\frac{1}{\varphi_{j}(t_{j})}
\sum\limits_{s_{m}=n_{m}+1}^{\infty}...\sum\limits_{s_{1}=n_{1}+1}^{\infty}
\|\delta_{\bar{s}}(f)\|_{X(\bar{\varphi})}^{*}+\right.
 $$
 $$
\left.+\sum\limits_{e\subset \{1,...,m\}}
\prod\limits_{j\notin e} \frac{1}{\varphi_{j}(t_{j})} \sum\limits_{\bar{s}\in
G_{e}(\bar{n})}\prod\limits_{j\in
e}\frac{1}{\varphi_{j}(2^{-n_{j}})}\|\delta_{\bar{s}}(f)
\|_{X(\bar{\varphi})}^{*}\right\},
 $$
for $\bar{t}=(t_{1},...,t_{m})\in (2^{-n_{1}-1},2^{-n_{1}}]\times...\times
 (2^{-n_{m}-1},2^{-n_{m}}],$  $n_{j}=1,2,...;$ \\ $j=1,...,m.$}
\end{theorem}
{\bf Proof.} $E_{j}\subset [0,2\pi)$ be a Lebesgue  measurable
subset. Then, by the properties of the integral we get
 $$
\int_{E_{m}}...\int_{E_{1}}|f(x_{1},...,x_{m})|dx_{1}...dx_{m}\leq
\int_{E_{m}}...\int_{E_{1}}|f(\bar x)-U_{\bar n}(f,\bar{x})|d{\bar x}+
 $$
 $$
+\int_{E_{m}}...\int_{E_{1}}|U_{\bar n}(f,\bar{x})|d\bar{x}. \eqno(6)
 $$
Using H\"{o}lder's integral inequality we obtain
 $$
\int_{E_{m}}...\int_{E_{1}}|f(\bar x)-U_{\bar n}(f,\bar{x})|d\bar{x}\leq
C\prod_{j=1}^{m}\frac{|E_{j}|}{\varphi_{j}(|E_{j}|)}
\|f-U_{\bar n}(f)\|_{X(\bar{\varphi})}^{*}.
\eqno(7)
 $$
 Let $\forall e \subset \{1,...,m\}.$ Then by applying Lemma 1 we obtain
 $$
\int_{E_{m}}...\int_{E_{1}}\left|\sum_{\bar{s}\in G_{e}(\bar
n)}\delta_{\bar s}(f,\bar{x})\right|d\bar{x} \leq
 $$
 $$
 C\prod_{j\in e}|E_{j}| \prod_{j\notin e}
\frac{|E_{j}|}{\varphi_{j}(|E_{j}|)}\sum_{\bar{s}\in G_{e}(\bar n)}\prod_{j\in
e}\frac{1}{\varphi_{j}(2^{-s_{j}})}\|\delta_{\bar s}\|_{X(\bar{\varphi})}^{*}.
 \eqno(8)
 $$
Further, taking into account that $|E_{j}| \geq t_{j}$ and the properties of the function $\varphi_{j}$ from inequalities (6)--(8) we have
 $$
\prod_{j=1}^{m}|E_{j}|^{-1}\int_{E_{m}}...\int_{E_{1}}|f(\bar
x)-U_{\bar n}(f,\bar{x})|d\bar{x}\leq
 $$
 $$
\leq C\cdot
\left\{\prod_{j=1}^{m}\frac{1}{\varphi_{j}(t_{j})}\sum_{s_{m}=n_{m}+1}^{\infty}
...\sum_{s_{1}=n_{1}+1}^{\infty}\|\delta_{\bar s}(f)
\|_{X(\bar{\varphi})}^{*}+\right.
 $$
 $$
\left.+\sum_{e\subset \{1,...,m\}}\prod_{j\notin e}\frac{1}{\varphi_{j}(t_{j})}
\sum_{\bar{s}\in G_{e}(\bar n)}\prod_{j\in e}
\frac{1}{\varphi_{j}(2^{-s_{j}})}\|\delta_{\bar s}(f)
\|_{X(\bar{\varphi})}^{*}\right\}.
 $$
 This implies the assertion of the theorem.

\begin{theorem}\label{th2} 
{\it Let $1<\alpha_{\psi_{j}}\leq \beta_{\psi_{j}}<
\alpha_{\varphi_{j}}\leq \beta_{\varphi_{j}}<2$, $1\leq\tau_{j}<+\infty$,
$j=1,...,m$. If  $f\in X(\bar{\varphi})$ and 
 $$
\left\{\prod_{j=1}^{m}\frac{\psi_{j}(2^{-s_{j}})}{\varphi_{j}
(2^{-s_{j}})}\|\delta_{\bar s}(f)\|_{X(\bar{\varphi})}^{*}
\right\}_{\bar{s}\in{\mathbb Z}_{+}^{m}}\in
l_{\bar{\tau}},
 $$
 then $f\in L_{\bar{\psi},\bar{\tau}}^{*}$ and  the following inequality holds
 $$
\|f\|_{\bar{\psi},\bar{\tau}}^{*}\leq C
\left\|\left\{\prod_{j=1}^{m}\frac{\psi_{j}(2^{-s_{j}})}{\varphi_{j}
(2^{-s_{j}})}\|\delta_{\bar s}(f)\|_{X(\bar{\varphi})}^{*}
\right\}_{\bar{s}\in{\mathbb Z}_{+}^{m}}
\right\|_{l_{\bar\tau}}.
 $$}
\end{theorem}
{\bf Proof.} According to Lemma 2 \cite{4}, the following inequality holds:
 $$
f^{*_{1},...,*_{m}}(t_{1},...,t_{m})\leq \bar{f}(t_{1},...,t_{m})\equiv
 $$
 $$
\equiv \sup_{|E_{m}|\geq t_{m}}\int_{E_{m}}dx_{m}...
\sup_{|E_{1}|\geq t_{1}}\int_{E_{1}}|f(x_{1},...,x_{m})|dx_{1}.
 $$
Therefore 
 $
\|f\|_{\bar{\psi},\bar{\tau}}^{*}\leq C\|\bar{f}\|_{
\bar{\psi},\bar{\tau}}^{*},
\,\, 1\leq\tau_{j}<+\infty, j=1,...,m.
 $
 Taking into account the relation 
 $$
\int\limits_{2^{-n-1}}^{2^{-n}}\psi_{j}(t)\frac{dt}{t}\asymp \psi_{j}(2^{-n}).
 \eqno(9)
 $$
 and using Theorem 1, we have
 $$
\|f\|_{\bar{\psi},\bar{\tau}}^{*}\leq
C\left[\left\{\sum_{n_{m}=0}^{\infty}\left(\frac{\psi_{m}(2^{-n_{m}})}{
\varphi_{m}(2^{-n_{m}})}\right)^{\tau_{m}}
\left[...\right.\right.\right.
 $$
 $$
...\left.\left.\left[\sum_{n_{1}=0}
^{\infty} \left(\frac{\psi_{1}(2^{-n_{1}})}{
\varphi_{1}(2^{-n_{1}})}\right)^{\tau_{1}} \left(
\sum_{s_{m}=n_{m}+1}^{\infty}...\sum_{s_{1}=n_{1}+1}^{\infty}
\|\delta_{\bar s}(f)\|_{X(\bar{\varphi})}^{*}\right)^{\tau_{1}}
\right]^{\frac{\tau_{2}}{\tau_{1}}}...\right]^{\frac{\tau_{m}}{
\tau_{m-1}}} \right\}^{\frac{1}{\tau_{m}}}+
 $$
 $$
+\sum\limits_{e\subset\{1,...,m\}}\left\{\sum_{n_{m}=0}^{\infty}\int\limits_{2^{-n_{m}-1}}^{
2^{-n_{m}}}
\psi_{m}^{\tau_{m}}(t_{m})t_{m}^{-1}\left[...\left[\sum_{n_{1}=0}
^{\infty}\int\limits_{2^{-n_{1}-1}}^{2^{-n_{1}}}\psi_{1}^{\tau_{1}}(t_{1})
\cdot t_{1}^{-1} \times \right.\right.\right.
 $$
 $$
\left.\times\left(\prod\limits_{j\notin
e}\frac{1}{\varphi_{j}(t_{j})}\sum\limits_{\bar{s}\in
G_{e}(\bar{n})}\prod\limits_{j\in
e}\frac{1}{\varphi_{j}(2^{-s_{j}})}\|\delta_{\bar{s}}(f)
\|_{X(\bar{\varphi})}^{*}\right)^{\tau_{1}}dt_{1}
\right]^{\frac{\tau_{2}}{\tau_{1}}}...
 $$
 $$
...\left.\left.\left.\right]^{\frac{\tau_{m}}{
\tau_{m-1}}} dt_{m}\right\}^{\frac{1}{\tau_{m}}}\right]=
 C\cdot\left[J_{1}+\sum\limits_{e\subset\{1,...,m\}}J_{e}\right]. \eqno (10)
 $$
According to Lemma 4 and Lemma 5, the numbers
$$
a_{s_{j}} = \frac{\psi(2^{-s_{j}})}{\varphi(2^{-s_{j}})}
$$
and $a_{s_{j}} = \psi(2^{-s_{j}})$, $j=1,\ldots,m$ satisfy the conditions $ a) $ and $ b) $ of Lemma 6, respectively.

Let $e=\{1,...,i\}, i\leq m$. Then using relation (9) and
successively applying the triangle inequality, assertions $ a) $ and $ b) $ of Lemma 6, we will have
 $$
J_{e}\leq \left\{\sum\limits_{n_{m}=0}^{\infty}
\left(\frac{\psi_{m}(2^{-n_{m}})}{
\varphi_{m}(2^{-n_{m}})}\right)^{\tau_{m}}  \left[...\right. \right.
 $$
 $$
...\left[\sum_{n_{i+1}=0}^{\infty}
\left(\frac{\psi_{i+1}(2^{-n_{i+1}})}{
\varphi_{i+1}(2^{-n_{i+1}})}\right)^{\tau_{i+1}}  \left[
\sum_{n_{i}=0}^{\infty}\psi_{i}^{\tau_{i}}(2^{-n_{i}})\left[...\left[
\sum_{n_{1}=0}^{\infty}\psi_{1}^{\tau_{1}}(2^{-n_{1}})\times
\right. \right. \right. \right.
 $$
 $$
\times \left(\sum_{s_{m}=n_{m}+1}^{\infty}...\sum_{s_{i+1}=n_{i+1}+1}^{\infty}
\sum_{s_{i}=1}^{n_{i}}...\sum_{s_{1}=1}^{n_{1}}\prod_{j=1}^{i}
\frac{1}{\varphi_{j}(2^{-s_{j}})}\times\right.
 $$
 $$
\left.\left.\left.\left.\left.
\times \|\delta_{\bar{s}}(f)\|_{X(\bar{\varphi})}^{*}
\right)^{\tau_{1}}
\right]^{\frac{\tau_{2}}{\tau_{1}}}...\biggr]^{\frac{\tau_{i}}
{\tau_{i-1}}} \right]^{\frac{\tau_{i+1}}{\tau_{i}}}...
\right]^{\frac{\tau_{m}}{\tau_{m-1}}}\right\}^{\frac{1}{\tau_{m}}} \leq
 $$
 $$
\leq C\left\{\sum\limits_{n_{m}=0}^{\infty}
\left(\frac{\psi_{m}(2^{-n_{m}})}{
\varphi_{m}(2^{-n_{m}})}\right)^{\tau_{m}}
\left[\sum\limits_{n_{m-1}=0}^{\infty}
\left(\frac{\psi_{m-1}(2^{-n_{m-1}})}{
\varphi_{m-1}(2^{-n_{m-1}})}\right)^{\tau_{m-1}}...\right. \right.
 $$
 $$
...\left.\left.\left[\sum\limits_{n_{1}=0}^{\infty}
\left(\frac{\psi_{1}(2^{-n_{1}})}{
\varphi_{1}(2^{-n_{1}})}\right)^{\tau_{1}}
\left(\|\delta_{\bar{s}}(f)\|_{X(\bar{\varphi})}^{*}\right)^{\tau_{1}}
\right]^{\frac{\tau_{2}}{\tau_{1}}}...\right]^{\frac{\tau_{m}}{
\tau_{m-1}}}
\right\}^{\frac{1}{\tau_{m}}}. \eqno(11)
 $$
Let $e \subset \{1,2,\ldots , m\},$ $e \neq \emptyset$ and $j_{0}=\min e$, $k_{0}=\max e$. 
 If $\{j_{0}, j_{0}+1,\ldots , k_{0}\}\cap \bar{e} = \emptyset$, where $\bar{e}$ is the complement of the set $e$, then
$$
\sum\limits_{\bar{s}\in
G_{e}(\bar{n})}\prod\limits_{j\in
e}\frac{1}{\varphi_{j}(2^{-s_{j}})}\|\delta_{\bar{s}}(f)
\|_{X(\bar{\varphi})}^{*} = 
$$  
$$
\sum_{s_{m}=n_{m}+1}^{\infty}...\sum_{s_{k_{0}+1}=n_{k_{0}+1}+1}^{\infty}
\sum_{s_{k_{0}}=1}^{n_{k_{0}}}...\sum_{s_{j_{0}}=1}^{n_{j_{0}}}
\sum_{s_{j_{0}-1}=n_{j_{0}-1}+1}^{\infty}...\sum_{s_{1}=n_{1}+1}^{\infty}\prod\limits_{j\in
e}\frac{1}{\varphi_{j}(2^{-s_{j}})}\|\delta_{\bar{s}}(f)
\|_{X(\bar{\varphi})}^{*}
$$
If $\{j_{0}, j_{0}+1,\ldots , k_{0}\}\cap \bar{e} = \{l_{0},\ldots , l_{1}\}$, then
$$
\sum\limits_{\bar{s}\in
G_{e}(\bar{n})}\prod\limits_{j\in
e}\frac{1}{\varphi_{j}(2^{-s_{j}})}\|\delta_{\bar{s}}(f)
\|_{X(\bar{\varphi})}^{*} = 
$$
$$
\sum_{s_{m}=n_{m}+1}^{\infty}...\sum_{s_{k_{0}+1}=n_{k_{0}+1}+1}^{\infty}
\sum_{s_{k_{0}}=1}^{n_{k_{0}}}...\sum_{s_{l_{1}-1}=0}^{n_{l_{1}-1}}
\sum_{s_{l_{1}}=n_{l_{1}}+1}^{\infty}...\sum_{s_{l_{0}}=n_{l_{0}}+1}^{\infty} 
\sum_{s_{l_{0}+1}=1}^{n_{l_{0}+1}}...\sum_{s_{j_{0}}=1}^{n_{j_{0}}}
$$
$$
\sum_{s_{j_{0}-1}=n_{j_{0}-1}+1}^{\infty}...\sum_{s_{1}=n_{1}+1}^{\infty}
\prod\limits_{j\in
e}\frac{1}{\varphi_{j}(2^{-s_{j}})}\|\delta_{\bar{s}}(f)
\|_{X(\bar{\varphi})}^{*}.
$$
Now, using these equalities and successively applying the assertions $ a) $, $ b) $ of Lemma 6, we obtain that inequality (11) holds for an arbitrary non-empty subset of $e$. Further, applying assertions $ a) $ of Lemma 6 we also obtain
 $$
J_{1}\leq C\cdot \left\{\sum\limits_{n_{m}=0}^{\infty}
\left(\frac{\psi_{m}(2^{-n_{m}})}{
\varphi_{m}(2^{-n_{m}})}\right)^{\tau_{m}}
\left[...\left[\sum\limits_{n_{1}=0}^{\infty}
\left(\frac{\psi_{1}(2^{-n_{1}})}{
\varphi_{1}(2^{-n_{1}})}\right)^{\tau_{1}}
\|\delta_{\bar{s}}(f)\|_{X(\bar{\varphi})}^{*}
\right]^{\frac{\tau_{2}}{\tau_{1}}}...\right]^{\frac{\tau_{m}}{
\tau_{m-1}}}
\right\}^{\frac{1}{\tau_{m}}}. \eqno(12)
 $$
Now the assertion of the theorem follows from inequalities (10)--(12).

\begin{rem} In the case $X(\bar\varphi)=L_{\bar{\varphi}, \bar{\tau}}^{*}(\mathbb{T}^{m})$ --- the generalized Lorentz space, Theorem 2 was announced in \cite{25} no proof.
\end{rem}

\begin{theorem}\label{th3} 
{\it Let  $1\leq\theta_{j}\leq +\infty$, $1\leq\tau_{j}<+\infty$, $j=1,...,m$ and functions $\varphi_{j}$, $\psi_{j}$ satisfies the conditions  $1<\alpha_{\psi_{j}}
\leq\beta_{\psi_{j}}<\alpha_{\varphi_{j}}\leq\beta_{\varphi_{j}}<2$, $j=1,...,m$.
If 
 $$
\left[\sum_{s_{j}=0}^{\infty}\left(\frac{\psi_{j}(2^{-s_{j}})}
{\varphi_{j}(2^{-s_{j}})}
2^{-s_{j}\tau_{j}}\right)^{\varepsilon_{j}}\right]^{\frac{1}{\varepsilon_{j}}}<+\infty,
\eqno(13)
 $$
where $\varepsilon_{j}=\tau_{j}\beta_{j}'$, $\beta_{j}'=
\frac{\beta_{j}}{\beta_{j}-1}$, $j=1,...,m$, if $\beta_{j}=
\frac{\theta_{j}}{\tau_{j}}>1$ and $\varepsilon_{j}=+\infty$, if
$\theta_{j}\leq\tau_{j}$, $j=1,...,m,$ then the embedding
$S_{X(\bar{\varphi}), \bar{\theta}}^{\bar r}B\subset
L_{\bar{\psi},\bar{\tau}}^{*}(I^{m})$ and
 $$
\|f\|_{\bar{\psi}, \bar{\tau}}^{*}\leq C
\|f\|_{S_{X(\bar{\varphi}), \bar{\theta}}^{\bar r}B}.
 $$}
\end{theorem}
{\bf Proof.} If $\tau_{j}<\theta_{j}$, $j=1,...,m$, then applying
H\"{o}lder's inequality with exponents $\beta_{j}=\frac{\theta_{j}}{\tau_{j}}$,
 $\frac{1}{\beta_{j}}+
\frac{1}{\beta_{j}'}=1$, $j=1,...,m$ we obtain
 $$
\sigma_{\bar{\tau}, \bar{\theta}, \bar{\varphi}}(f)=\left\|\left\{
\prod_{j=1}^{m}\frac{\psi_{j}(2^{-s_{j}}
)}{\varphi_{j}(2^{-s_{j}})}\|\delta_{\bar s}(f)\|_{X(\bar{\varphi})}^{*}
\right\}_{\bar s \in \mathbb{Z}_{+}^{m}}\right\|_{l_{\bar\tau}}\leq
 $$
 $$
\leq \left\|\left\{
\prod_{j=1}^{m}\frac{\psi_{j}(2^{-s_{j}})}{\varphi_{j}(2^{-s_{j}})}2^{-s_{j}r_{j}}
\right\}_{\bar s \in \mathbb{Z}_{+}^{m}}\right\|_{l_{\bar\varepsilon}}
\left\|\left\{\prod_{j=1}^{m}2^{s_{j}r_{j}}\|\delta_{\bar s}(f)\|_{X(\bar{\varphi})}^{*}
\right\}_{\bar s \in \mathbb{Z}_{+}^{m}}\right\|_{l_{\bar\theta}}, \eqno(14)
 $$
where $\bar\varepsilon=(\varepsilon_{1},...,\varepsilon_{m}),$ $\varepsilon_{j}=\tau_{j}\beta_{j}',$  $j=1,...,m.$

If $\theta_{j}\leq\tau_{j},$ $j=1,...,m,$ then by Jensen's inequality will have
 $$
\sigma_{\bar{\tau},\bar{\theta}, \bar{\varphi}}(f) \leq
\left\|\left\{
\prod_{j=1}^{m}
2^{-s_{j}r_{j}}
\|\delta_{\bar s}(f)\|_{X(\bar{\varphi})}^{*}
\right\}_{\bar s \in \mathbb{Z}_{+}^{m}}\right\|_{l_{\bar\theta}}
\prod_{j=1}^{m}\sup_{s_{j}\in\mathbb{Z}_{+}}
\frac{\psi_{j}(2^{-s_{j}})}
{\varphi_{j}(2^{-s_{j}})}2^{-s_{j}r_{j}}. \eqno(15)
 $$
By conditions (13) and Theorem 2, inequalities (14) and (15) imply the assertion
theorems.

\begin{theorem}\label{th4}  
 {\it Let functions  $\psi_{j}$ satisfies the conditions  $1< 2^{1/\lambda_{j}} < \alpha_{\psi_{j}}\leq \beta_{\psi_{j}}<2$, $1<\tau_{j}<+\infty$, $j=1,...,m.$
If $f\in L_{\bar{\psi}, \bar{\tau}}^{*}(\mathbb{T}^{m})$ and
 $$
f(\bar x)\sim \sum_{\bar{s}\in\mathbb{Z}_{+}^{m}}b_{\bar s}
\sum_{\bar{k}\in\rho(\bar s)}e^{i\langle\bar{k},\bar{x}\rangle},
 $$
 then the inequality holds
 $$
\|f\|_{\bar{\psi},\bar{\tau}}^{*}\geq C\times
 $$
 $$
\times 
\left\{\sum_{s_{m}=1}^{\infty}2^{s_{m}\frac{\tau_{m}}{\lambda_{m}}}\psi_{m}^{\tau_{m}}(2^{-s_{m}})\left[...\left[
\sum_{s_{1}=1}^{\infty}2^{s_{1}\frac{\tau_{1}}{\lambda_{1}}}\psi_{1}^{\tau_{1}}(2^{-s_{1}})\left(\|\delta_{\bar
s}(f)\|_{\bar{\lambda}, \bar{\tau}}^{*}
\right)^{\tau_{1}}\right]^{\frac{\tau_{2}}{\tau_{1}}}...
\right]^{\frac{\tau_{m}}{\tau_{m-1}}}\right\}^{\frac{1}{\tau_{m}}},
 $$
where  $b_{\bar s}$---real numbers.
 }
\end{theorem}

{\bf Proof.} Let $S_{2^{\nu},...,2^{\nu}}(f,\bar{x})$ be rectangular
partial sum of the Fourier series of a function 
$f\in L_{\bar{\psi}, \bar{\tau}}^{*}(I^{m}).$ It is known that (see \cite{4}) 
  $$
\|f\|_{\bar{\psi}, \bar{\theta}}^{*} \asymp \sup_{{}_{\|g\|_{\bar{\tilde{\psi}},\bar{\tau}'}^{*}\leq 1}^{g\in L_{\bar{\tilde{\psi}}',\bar{\tau}'}^{*}}} \int_{I^{m}}
f(2\pi\bar{x})g(2\pi\bar{x})d\bar{x}, \eqno(16)
 $$
where  $\bar{\tau}'= (\tau_{1}',...,\tau_{m}'),$ 
$\frac{1}{\tau_{j}}+\frac{1}{\tau_{j}'}=1$, $j=1,...,m$ and $\bar{\tilde{\psi}}(t)=(\tilde{\psi}_{1}(t),...,\tilde{\psi}_{m}(t))$, 
 $
 \tilde{\psi}_{j}(t)=\frac{t}{\psi_{j}(t)}, 
 $
 for $t\in (0, 1]$ and $\tilde{\psi}_{j}(0)=0$, $j=1,...,m$. 
We will introduce the notation 
 $$
\sigma_{\nu}(f)_{\bar{\tau}_{j}}=
\left\{\sum_{s_{j}=1}^{\nu-1}2^{\frac{s_{j}\tau_{j}}{\lambda_{j}}}\psi_{j}^{\tau_{j}}(2^{-s_{j}})\left[...\left[
\sum_{s_{1}=1}^{\nu-1}2^{\frac{s_{1}\tau_{1}}{\lambda_{1}}}\psi_{1}^{\tau_{1}}(2^{-s_{1}})\left(\|\delta_{\bar s}(f)\|_{\bar{\lambda}, \bar{\tau}}^{*}
\right)^{\tau_{1}}\right]^{\frac{\tau_{2}}{\tau_{1}}}...
\right]^{\frac{\tau_{j}}{\tau_{j-1}}}\right\}^{\frac{1}{\tau_{j}}}.
 $$
We consider trigonometric polynomial
 $$
g_{\nu}(\bar x)=\sum_{s_{m}=1}^{\nu-1}...\sum_{s_{1}=1}^{\nu-1}b_{\bar{s},\nu}
\sum_{\bar{k}\in\rho(\bar s)}e^{i\langle\bar{k},\bar{x}\rangle},
 $$
where 
 $$
b_{\bar{s},\nu}=\left\|\left\{\prod_{j=1}^{m}2^{\frac{s_{j}}{\lambda_{j}}}\psi_{j}(2^{-s_{j}})\|\delta_{\bar s}(f)\|_{\bar{\lambda}, \bar{\tau}}^{*}
\right\}_{\bar{s}=\bar{1}}^{\bar{\nu}}\right\|_{l_{\bar\tau}}^{
-\frac{\tau_{m}}{\tau'_{m}}} \prod_{j=1}^{m-1}(\sigma_{\nu}
(f)_{\bar{\tau}_{j}})^{\tau_{j+1}-\tau_{j}}
 $$
 $$
\times \prod_{j=1}^{m}(2^{s_{j}}\psi_{j}(2^{-s_{j}}))^{\tau_{j}}
\left(\prod_{j=1}^{m}2^{s_{j}}\right)^{-1}
\prod_{j=2}^{m}2^{s_{j}(1-\frac{1}{\lambda_{j}})(\tau_{1}-\tau_{j})}
|b_{\bar s}(f)|^{\tau_{1}-1}sign(b_{\bar s}(f))
 $$
and $\bar{\tau}_{j}=(\tau_{1}, \ldots , \tau_{j})$.
Then, taking into account the orthogonality of the trigonometric system, we have
$$ 
 \int_{I^{m}}
f(2\pi\bar{x})g(2\pi\bar{x})d\bar{x} = \int_{I^{m}}
S_{2^{\nu},...,2^{\nu}}(f,2\pi\bar{x}) g_{\nu}(2\pi\bar x)d\bar{x} 
$$
$$
=
\sum_{s_{m}=1}^{\nu-1}...\sum_{s_{1}=1}^{\nu-1}
\int_{I^{m}}\delta_{\bar s}(f, 2\pi\bar{x})g_{\nu}(2\pi\bar x)d\bar{x}.
 \eqno(17)
 $$
 We will prove that 
 $\|g_{\nu}\|_{\bar{\tilde{\psi}}, \bar{\tau}'}^{*} \leq C_{0},$ 
 where $C_{0}$ is some positive constant independent of $\nu$.
 Taking into account the relation (see (5))
 $$
\Bigl
\|\sum_{\bar{k}\in\rho(\bar s)}e^{i\langle\bar{k},\bar{x}\rangle}
\Bigr\|_{\bar{\psi},\bar{\tau}}^{*}\asymp \prod_{j=1}^{m}
2^{s_{j}}\psi_{j}(2^{-s_{j}}),\,\,  1< \tau_{j}<+\infty, \,\, 1< \alpha_{\psi_{j}} \leq  \beta_{\psi_{j}}<2,     \eqno(18)
 $$
 it is easy to verify the following inequality
 $$
\left(\|\delta_{\bar s}(g_{\nu})\|_{\bar{\lambda}', \bar{\tau}'}^{*}
\right)^{\tau'_{1}}=\Bigl(|b_{s,\nu}| \Bigl\|
\sum_{\bar{k}\in\rho(\bar s)}e^{i\langle\bar{k},\bar{x}\rangle}
\Bigr\|_{\bar{\lambda}', \bar{\tau}'}^{*}\Bigr)^{\tau'_{1}}\leq
 $$
 $$
\leq C\prod_{j=2}^{m}\Bigl(2^{\frac{s_{j}}{\lambda_{j}}}\psi_{j}(2^{-s_{j}})\Bigr)^{\tau_{j}\tau'_{1}}
\left(2^{\frac{s_{1}}{\lambda_{1}}}\psi_{1}(2^{-s_{1}})\right)^{(\tau_{1}+\tau'_{1})}\left(\|\delta_{\bar s}(f)\|_{\bar{\lambda}, \bar{\tau}}^{*}\right)^{\tau_{1}}
 $$
 $$
\times \left(\prod_{j=1}^{m}(\sigma_{\nu}
(f)_{\bar{\tau}_{j}})^{\tau_{j+1}-\tau_{j}}
 \left\|\left\{\prod_{j=1}^{m}2^{\frac{s_{j}}{\lambda_{j}}}\psi_{j}(2^{-s_{j}})\|\delta_{\bar s}(f)\|_{\bar{\lambda}, \bar{\tau}}^{*}
\right\}_{\bar{s}=\bar{1}}^{\bar{\nu}}\right\|_{l_{\bar\tau}}^{
-\frac{\tau_{m}}{\tau'_{m}}}\right)^{\tau'_{1}}, \eqno(19)
 $$
 where $\bar{\nu} = (\nu, \ldots , \nu )$.
Further, by using inequality (19), we obtain
 $$
J(g_{\nu}):=
\left\{\sum_{s_{m}=1}^{\nu-1}(2^{\frac{s_{m}}{\lambda_{m}^{'}}}\tilde{\psi}_{m}(2^{-s_{m}}))^{\tau'_{m}}\left[...\left[
\sum_{s_{1}=1}^{\nu-1}(2^{\frac{s_{1}}{\lambda_{1}^{'}}}\tilde{\psi}_{1}(2^{-s_{1}}))^{\tau'_{1}}\left(\|\delta_{\bar s}(f)\|
_{\bar{\lambda}', \bar{\tau}'}^{*}
\right)^{\tau'_{1}}\right]^{\frac{\tau'_{2}}{\tau'_{1}}}...
\right]^{\frac{\tau'_{m}}{\tau'_{m-1}}}\right\}^{\frac{1}{\tau'_{m}}}
 $$
 $$
\leq C \left\|\left\{\prod_{j=1}^{m}2^{\frac{s_{j}}{\lambda_{j}}}\psi_{j}(2^{-s_{j}})\|\delta_{\bar s}(f)\|_{\bar{\lambda}, \bar{\tau}}^{*}
\right\}_{\bar{s}\in\mathbb{Z}_{+}^{m}}\right\|_{l_{\bar\tau}}^{
-\frac{\tau_{m}}{\tau'_{m}}}
 $$
 $$
\times
\left\{\sum_{s_{m}=1}^{\nu-1}(2^{\frac{s_{m}}{\lambda_{m}^{'}}}\tilde{\psi}_{m}(2^{-s_{m}}))^{\tau'_{m}}\left[...\left[
\sum_{s_{1}=1}^{\nu-1}(2^{\frac{s_{1}}{\lambda_{1}^{'}}}\tilde{\psi}_{1}(2^{-s_{1}}))^{\tau'_{1}}\prod_{j=2}^{m}
(2^{\frac{s_{j}}{\lambda_{j}}}\psi_{j}(2^{-s_{j}}))^{\tau_{j}\tau'_{1}}
\right.\right.\right.
 $$
 $$
\left.\left.\left.\times (2^{\frac{s_{1}}{\lambda_{1}}}\psi_{1}(2^{-s_{1}}))^{(\tau_{1}+\tau'_{1})}
\left(\|\delta_{\bar s}(f)\|_{\bar{\lambda}, \bar{\tau}}^{*}\right)^{\tau_{1}}
\left(\prod_{j=1}^{m-1}(\sigma_{\nu}
(f)_{\bar{\tau}_{j}})^{\tau_{j+1}-\tau_{j}}\right)^{\tau'_{1}}
\right]^{\frac{\tau'_{2}}{\tau'_{1}}}...
\right]^{\frac{\tau'_{m}}{\tau'_{m-1}}}\right\}^{\frac{1}{\tau'_{m}}}.
\eqno(20)
 $$
Since 
   $\frac{1}{\lambda_{j}^{'}}=1-\frac{1}{\lambda_{j}}$, $\tau_{j}\tau'_{j}=\tau_{1}+\tau'_{1}$ and $\tilde{\psi}_{j}(t)=\frac{t}{\psi_{j}(t)}$, $j=1, \ldots, m$ then 
 $$
 \Bigl(2^{\frac{s_{j}}{\lambda_{j}^{'}}}\tilde{\psi}_{j}(2^{-s_{j}})\Bigr)^{\tau'_{j}} \Bigl(2^{\frac{s_{j}}{\lambda_{j}}}\psi_{j}(2^{-s_{j}})\Bigr)^{\tau_{j}\tau'_{j}} = \Bigl(\frac{2^{-\frac{s_{j}}{\lambda_{j}}}}{\psi_{j}(2^{-s_{j}})}\Bigr)^{\tau'_{j}}   \Bigl(2^{\frac{s_{j}}{\lambda_{j}}}\psi_{j}(2^{-s_{j}})\Bigr)^{\tau_{j}+\tau'_{j}} 
 $$ 
 $$
 = \Bigl(2^{\frac{s_{j}}{\lambda_{j}}}\psi_{j}(2^{-s_{j}})\Bigr)^{\tau_{j}}, \,\, j=1, \ldots, m.
 $$
 Now, considering this equality and the definition of numbers
$\sigma_{\nu}(f)_{\bar{\tau}_{j}}$ inequality (20) continues
 $$
J(g_{\nu})\leq C
\left\|\left\{\prod_{j=1}^{m}2^{\frac{s_{j}}{\lambda_{j}}}\psi_{j}(2^{-s_{j}})\|\delta_{\bar s}(f)\|_{\bar{\lambda}, \bar{\tau}}^{*}
\right\}_{\bar{s}\in\mathbb{Z}_{+}^{m}}\right\|_{l_{\bar\tau}}^{
-\frac{\tau_{m}}{\tau'_{m}}}
 $$
 $$
\times \left\{\sum_{s_{m}=1}^{\nu-1}\Bigl(2^{\frac{s_{m}}{\lambda_{m}}}\psi_{m}(2^{-s_{m}})\Bigr)^{\tau_{m}}\left[...\left[
\sum_{s_{1}=1}^{\nu-1}\Bigl(2^{\frac{s_{1}}{\lambda_{1}}}\psi_{1}(2^{-s_{1}})\Bigr)^{\tau_{1}}
\left(\|\delta_{\bar s}(f)\|_{\bar{\lambda}, \bar{\tau}}^{*}\right)^{\tau_{1}}
\times\right.\right.\right.
 $$
 $$
\left.\left.\left.\times
\left(\prod_{j=1}^{m-1}(\sigma_{\nu}
(f)_{\bar{\tau}_{j}})^{\tau_{j+1}-\tau_{j}}\right)^{\tau'_{1}}
\right]^{\frac{\tau'_{2}}{\tau'_{1}}}...
\right]^{\frac{\tau'_{m}}{\tau'_{m-1}}}\right\}^{\frac{1}{\tau'_{m}}}=
 $$
 $$
=
 \left\{\sum_{s_{m}=1}^{\nu-1}\Bigl(2^{\frac{s_{m}}{\lambda_{m}}}\psi_{m}(2^{-s_{m}})\Bigr)^{\tau_{m}}\left[...\left[
\sum_{s_{1}=1}^{\nu-1}\Bigl(2^{\frac{s_{1}}{\lambda_{1}}}\psi_{1}(2^{-s_{1}})\Bigr)^{\tau_{1}}
\left(\|\delta_{\bar s}(f)\|_{\bar{\lambda}, \bar{\tau}}^{*}\right)^{\tau_{1}}
\right]^{\frac{\tau_{2}}{\tau_{1}}}...
\right]^{\frac{\tau_{m}}{\tau_{m-1}}}\right\}^{\frac{1}{\tau_{m}}}
 $$
 $$
 \times C
\left\|\left\{\prod_{j=1}^{m}2^{\frac{s_{j}}{\lambda_{j}}}\psi_{j}(2^{-s_{j}})\|\delta_{\bar s}(f)\|_{\bar{\lambda}, \bar{\tau}}^{*}
\right\}_{\bar{s}\in\mathbb{Z}_{+}^{m}}\right\|_{l_{\bar\tau}}^{
-\frac{\tau_{m}}{\tau'_{m}}}
 \leq C. \eqno(21)
 $$
Since by the hypothesis of the theorem 
 $1< 2^{\frac{1}{\lambda_{j}}} < \alpha_{\psi_{j}}$,  then 
 $\beta_{\tilde{\psi_{j}}} < 2^{\frac{1}{\lambda_{j}^{'}}}$,  $j=1,...,m.$
Therefore, according to Theorem 2 and inequality (21), we obtain
 $$
\|g_{\nu}\|_{\bar{\tilde{\psi}}, \bar{\tau}'}^{*}\leq C J(g_{\nu})
\leq C_{0}.
 $$
 Hence the function  
 $\varphi_{\nu}=\frac{1}{C_{0}}g_{\nu}
\in L_{\bar{\tilde{\psi}}, \bar{\tau}'}^{*}(I^{m})$ and
$\|\varphi_{\nu}\|_{\bar{\tilde{\psi}}, \bar{\tau}'}^{*}\leq 1.$
Further, according to the orthogonality property, we have
 $$
\int_{I^{m}}\delta_{\bar s}(f,\bar{x})g_{\nu}(\bar x)d\bar{x}=
\sum_{l_{m}=1}^{\nu-1}...\sum_{l_{1}=1}^{\nu-1}b_{\bar s}(f)
b_{\bar{l},\nu}\times
 $$
 $$
\times \int_{I^{m}}\left|\sum_{\bar{k}\in\rho(\bar s)}
e^{i\langle\bar{k},\bar{x}\rangle}\right|^{2}d\bar{x}=(2\pi)^{d}\cdot
b_{\bar s}(f)b_{\bar{s},\nu}\cdot \prod_{j=1}^{m}2^{s_{j}-1}. \eqno(22)
 $$
From the definition of the numbers $b_{\bar{s},\nu}$ it follows that
 $$
\prod_{j=1}^{m}2^{s_{j}}b_{\bar s}(f)b_{\bar{s},\nu}=\prod_{j=1}^{m}2^{s_{j}}
\left\|\left\{\prod_{j=1}^{m}2^{\frac{s_{j}}{\lambda_{j}}}\psi_{j}(2^{-s_{j}})\|\delta_{\bar s}(f)\|_{\bar{\lambda}, \bar{\tau}}^{*}
\right\}_{\bar{s}=\bar{1}}^{\bar{\nu}}\right\|_{l_{\bar\tau}}^{
-\frac{\tau_{m}}{\tau'_{m}}} \prod_{j=1}^{m-1}(\sigma_{\nu}
(f)_{\bar{\tau}_{j}})^{\tau_{j+1}-\tau_{j}}
 $$
 $$
\times \prod_{j=1}^{m}2^{s_{j}\theta_{j}}\psi_{j}^{\tau_{j}}(2^{-s_{j}})
\left(\prod_{j=1}^{m}2^{s_{j}}\right)^{-1}
\prod_{j=2}^{m}2^{s_{j}(1-\frac{1}{\lambda_{j}})(\tau_{1}-\tau_{j})}
|b_{\bar s}(f)|^{\tau_{1}}=
 $$
 $$
=\prod_{j=1}^{m-1}(\sigma_{\nu}
(f)_{\bar{\tau}_{j}})^{\tau_{j+1}-\tau_{j}}
\left\|\left\{\prod_{j=1}^{m}2^{\frac{s_{j}}{\lambda_{j}}}\psi_{j}(2^{-s_{j}})\|\delta_{\bar s}(f)\|_{\bar{\lambda}, \bar{\tau}}^{*}
\right\}_{\bar{s}=\bar{1}}^{\bar{\nu}}\right\|_{l_{\bar\tau}}^{
-\frac{\tau_{m}}{\tau'_{m}}}
\prod_{j=1}^{m}2^{s_{j}\tau_{j}}\psi_{j}^{\tau_{j}}(2^{-s_{j}})
 $$
 $$
\times \prod_{j=2}^{m}2^{-s_{j}(1-\frac{1}{\lambda_{j}})\tau_{j}}
2^{-s_{1}(1-\frac{1}{\lambda_{1}})\tau_{1}} \left(\prod_{j=1}^{m}
2^{s_{j}(1-\frac{1}{\lambda_{j}})}|b_{\bar s}(f)|\right)^{\tau_{1}}\geq
 $$
 $$
\geq C(\tau, \lambda)
\left\|\left\{\prod_{j=1}^{m}2^{\frac{s_{j}}{\lambda_{j}}}\psi_{j}(2^{-s_{j}})\|\delta_{\bar s}(f)\|_{\bar{\lambda}, \bar{\tau}}^{*}
\right\}_{\bar{s}=\bar{1}}^{\bar{\nu}}\right\|_{l_{\bar\tau}}^{
-\frac{\tau_{m}}{\tau'_{m}}}
\prod_{j=1}^{m-1}(\sigma_{\nu}
(f)_{\bar{\tau}_{j}})^{\tau_{j+1}-\tau_{j}}
 $$
 $$
\times \prod_{j=1}^{m}\Bigl(2^{\frac{s_{j}}{\lambda_{j}}}\psi_{j}(2^{-s_{j}})\Bigr)^{\tau_{j}}
\left(\|\delta_{\bar s}(f)\|_{\bar{\lambda}, \bar{\tau}}^{*}\right)^{\tau_{1}}.
\eqno(23)
 $$
Now, taking into account (22) and (23), we obtain
 $$
\sup_{{}_{\|g\|_{\bar{\tilde{\psi}}, \bar{\tau}'}^{*}
\leq 1}^{g\in L_{\bar{\tilde{\psi}}, \bar{\tau}'}^{*}}}
\sum_{s_{m}=1}^{\nu-1}...\sum_{s_{1}=1}^{\nu-1}
\int_{I^{m}}\delta_{\bar s}(f,\bar{x}) g(\bar x)d\bar{x}\geq
\sum_{s_{m}=1}^{\nu-1}...\sum_{s_{1}=1}^{\nu-1}
\int_{I^{m}}\delta_{\bar s}(f,\bar{x})\varphi_{\nu}(\bar x) d\bar{x}=
 $$
 $$
=C
\sum_{s_{m}=1}^{\nu-1}...\sum_{s_{1}=1}^{\nu-1}
\int_{I^{m}}\delta_{\bar s}(f,\bar{x}) g(\bar x)d\bar{x}\geq
\sum_{s_{m}=1}^{\nu-1}...\sum_{s_{1}=1}^{\nu-1}
\int_{I^{m}}\delta_{\bar s}(f,\bar{x})g_{\nu}(\bar x) d\bar{x} 
 $$
 $$
\geq C
\left\|\left\{\prod_{j=1}^{m}2^{\frac{s_{j}}{\tau_{j}}}\psi_{j}(2^{-s_{j}})\|\delta_{\bar s}(f)\|_{\bar{\tau},\bar{\theta}}^{*}
\right\}_{\bar{s}=\bar{1}}^{\bar{\nu}}\right\|_{l_{\bar\theta}}^{
-\frac{\theta_{m}}{\theta'_{m}}}
\prod_{j=1}^{m-1}(\sigma_{\nu}
(f)_{\bar{\theta}_{j}})^{\theta_{j+1}-\theta_{j}}
 $$
 $$
\times \sum_{s_{m}=1}^{\nu-1}...\sum_{s_{1}=1}^{\nu-1}
\prod_{j=1}^{m}\Bigl(2^{\frac{s_{j}}{\tau_{j}}}\psi_{j}(2^{-s_{j}})\Bigr)^{\theta_{j}}
\left(\|\delta_{\bar s}(f)\|_{\bar{\lambda}, \bar{\tau}}^{*}\right)^{\tau_{1}}=
 $$
 $$
=C
\left\|\left\{\prod_{j=1}^{m}2^{\frac{s_{j}}{\lambda_{j}}}\psi_{j}(2^{-s_{j}})\|\delta_{\bar s}(f)\|_{\bar{\lambda}, \bar{\tau}}^{*}
\right\}_{\bar{s}=\bar{1}}^{\bar{\nu}}\right\|_{l_{\bar\tau}}. \eqno(24)
 $$

It follows from inequalities (17) and (24) that
  $$
\|f\|_{\bar{\psi}, \bar{\tau}}^{*}\geq C
\left\|\left\{\prod_{j=1}^{m}2^{\frac{s_{j}}{\lambda_{j}}}\psi_{j}(2^{-s_{j}})\|\delta_{\bar s}(f)\|_{\bar{\lambda}, \bar{\tau}}^{*}
\right\}_{\bar{s}\in\mathbb{Z}_{+}^{m}}\right\|_{l_{\bar\tau}}.
  $$
The theorem is proved. 

\begin{rem}  In the case 
 $\psi_{j}(t)=t^{\frac{1}{q_{j}}}$, $\varphi_{j}(t)=t^{\frac{1}{p_{j}}}$, $j=1,\ldots , m$ Theorem 2 and Theorem 4 were proved in \cite{14}. 
\end{rem}

\begin{theorem}\label{th5} 
{\it Let $1\leq\theta_{j}<+\infty,$ $1\leq\tau_{j}<+\infty,$ $j=1,...,m$ and 
functions $\varphi_{j},$ $\psi_{j}$ satisfies the conditions $1<\alpha_{\psi_{j}}
\leq\beta_{\psi_{j}}<\alpha_{\varphi_{j}}\leq\beta_{\varphi_{j}}<2,$ $j=1,...,m$
and (13).

1) If $1\leq \tau_{j}<\theta_{j}<+\infty,$ $j=1,...,m,$ then
 $$
E_{n}^{\bar \gamma}(S_{X(\bar{\varphi}), \bar\theta}^{\bar r}B)
_{\bar{\psi}, \bar{\tau}}\leq
C\left\|\left\{\prod_{j=1}^{m}2^{-s_{j}r_{j}}
\mu_{j}(s_{j})\right\}_{\bar{s}\in Y^{m}(\bar{\gamma},n)}
\right\|_{\bar\varepsilon},
 $$
where $\bar{\varepsilon}=(\varepsilon_{1},...,\varepsilon_{m}),$ 
$\varepsilon_{j}=\tau_{j}\beta_{j}',$  $\frac{1}{\beta_{j}}+
\frac{1}{\beta_{j}'}=1,$  $\beta_{j}=\frac{\theta_{j}}{\tau_{j}}.$

2) If $1\leq \theta_{j}\leq \tau_{j}<+\infty,$  $j=1,...,m,$ then
 $$
E_{n}^{\bar \gamma}(S_{X(\bar{\varphi}), \bar\theta}^{\bar r}B)
_{\bar{\psi},\bar{\tau}}\leq
C\sup\left\{\prod_{j=1}^{m}2^{-s_{j}r_{j}}
\mu_{j}(s_{j}): \quad \bar{s}\in\mathbb{Z}_{+}^{m},
\quad \langle\bar{s},\bar{\gamma}\rangle\geq n\right\}.
 $$
}
\end{theorem}
{\bf Proof.} Let $f\in S_{X(\bar{\varphi}), \bar\theta}^{\bar r}B.$ Then by  Theorem 3 we obtain  
 $$
\|f-S_{n}^{(\bar\gamma)}(f)\|_{\bar{\psi}, \bar{\tau}}^{*}\leq
C\left\|\left\{\prod_{j=1}^{m}
\mu_{j}(s_{j})\|\delta_{\bar s}(f-S_{n}^{(\bar\gamma)}(f))\|_{X(\bar{\varphi})}^{*}\right\}_{\bar s \in \mathbb{Z}_{+}^{m}}\right\|_{l_{\bar\tau}}. \eqno(25)
 $$
Since 
 $\delta_{\bar s}(f-S_{n}^{(\bar\gamma)}(f))=0$, $\bar{s}
\notin Y^{m}(\bar{\gamma},n)$ and
$\delta_{\bar s}(f-S_{n}^{(\bar\gamma)}(f))=\delta_{\bar s}(f)$, $\bar{s}
\in Y^{m}(\bar{\gamma},n),$ then from (25) we obtain 
 $$
\|f-S_{n}^{(\bar\gamma)}(f)\|_{\bar{\psi}, \bar{\tau}}^{*}\leq
C(\theta, m)\left\|\left\{\prod_{j=1}^{m}\mu_{j}(s_{j})
\|\delta_{\bar s}(f)\|_{X(\bar{\varphi})}^{*}
\right\}_{\bar{s}\in Y^{m}(\bar{\gamma},n)}\right\|_{l_{\bar\tau}}. \eqno(26)
 $$
We put 
 $$
b_{\bar s}(n)=\prod_{j=1}^{m}\frac{\psi_{j}\left(2^{-s_{j}}\right)}
{\varphi_{j}\left(2^{-s_{j}}\right)} , \,\, \text{????} \quad \bar{s}\in Y^{m}(\bar{\gamma},n),
 $$
$b_{\bar s}(n)=0$, $\bar{s}\notin Y^{m}(\bar{\gamma},n)$.

We will prove item 1). Since $1\leq \tau_{j}<\theta_{j}<+\infty,$ $j = 1, ...,m,$ then applying H\"{o}lder's inequality with exponents
$\beta_{j}=\frac{\tau_{j}}{\theta_{j}}$,
 $\frac{1}{\beta_{j}}+
\frac{1}{\beta_{j}'}=1$, $j=1,...,m$ we get 
 $$
\sigma_{1}(f,n)=\left\|\left\{
\prod_{j=1}^{m}
\mu_{j}(s_{j})
\|\delta_{\bar s}(f)\|_{X(\bar{\varphi})}^{*}
\right\}_{\bar{s}\in Y^{m}(\bar{\gamma},n)}\right\|_{
l_{\bar\tau}}\leq C \left\|\left\{b_{\bar s}(n)
\right\}_{\bar{s}\in\mathbb{Z}_{+}^{m}}\right\|_{l_{\bar\varepsilon}}, 
\eqno(27)
 $$
where $\bar{\varepsilon}=(\varepsilon_{1},...,\varepsilon_{m}),$
 $\varepsilon_{j}=\tau_{j}\beta_{j}',$ \,\, $j=1,...,m.$

Let us prove item 2). If $\theta_{j} \leq \tau_{j}<+\infty,$ $j=1,..., m,$ then according to Jensen's inequality (see \cite[Ch. 3, Sec. 3]{6}) we obtain
 $$
\sigma_{1}(f,n)\leq \|f\|_{S_{X(\bar{\varphi}), \bar{\theta}}^{\bar r}B}
 \sup_{\bar{s}\in Y^{m}(\bar{\gamma},n)}
\prod_{j=1}^{m}\mu_{j}(s_{j})2^{-s_{j}r_{j}}. \eqno(28)
 $$
Inequalities (26)--(28) imply the statements of items 1) and 2) of Theorem 5.

\begin{rem} In the case $\varphi_{j}(t)=t^{1/p}$, $p_{j}=\tau_{j}^{(1)}=p $ and $ \psi_{j}(t)=t^{1/q}$, $q_{j}=\tau_{j}^{(2)}=q$, $1\leq \theta_{j}=\theta\leq \infty$ for $j = 1, ..., m$ Theorem 5 implies the previously known results of Ya.S. Bugrova, E.M. Galeeva, V.N. Temlyakov and A.S. Romanyuk (see, for example, the bibliography in \cite{11}, \cite{12}, \cite{22}).

Further, From Theorem 5 with $\varphi_{j}(t)=t^{1/p_{j}}$, $\psi_{j}(t)=t^{1/q_{j}}$, $1< p_{j}, q_{j}< \infty$, $ j = 1, ..., m $ and $\gamma_{j}^{'} =\gamma_{j}=1$ for $j = 1, ..., \nu $ and $\gamma_{j}^{'} < \gamma_{j}$  $j=\nu+1,...,m$ follow Theorem 2 in \cite{13} ( also see \cite{14}) and for  $1\leq \gamma_{j}^{'} \leq \gamma_{j}$ for $ j = 1, ..., m $, Theorem 1 in \cite{16} ( see also \cite{18}).
\end{rem}

\begin{rem}
In the case $X(\bar\varphi)=L_{\bar{\varphi}, \bar{\tau}}^{*}(\mathbb{T}^{m})$ --- the generalized Lorentz space, Theorem 4 and Theorem 5 proved in \cite{25}.
\end{rem}

This work was supported by a grant
  Ministry of Education and Science of the Republic of Kazakhstan (Project AP 08855579).

\end{document}